\documentclass[11pt]{article}
\usepackage{amsmath,amsfonts,amssymb,latexsym}

\newtheorem{theorem}{Theorem}[section]
\newtheorem{proposition}[theorem]{Proposition}

\newtheorem{coro}[theorem]{Corollary}
\newtheorem{remark}[theorem]{Remark}
\newtheorem{ex}[theorem]{Example}

\def\R{{\mathbb R}}
\def\N{{\mathbb N}}
\newcommand{\Scal}{\mathcal{S}}
\newcommand{\Hcal}{\mathcal{H}}
\newcommand{\D}{\mathrm{D}}
\newcommand{\xR}{{]}{-\infty},+\infty]}
\newcommand{\Rex}{\xR}

\newcommand{\ps}{\smallbreak}
\newcommand{\proof}{\noindent{\bf Proof. }}
\newcommand{\cqfd}{\mbox{}\nolinebreak\hfill\rule{2mm}{2mm}\medbreak\par}

\newcommand{\esum}{+_{\rm e}}
\newcommand{\dom} {{\rm dom} \kern.15em}
\newcommand{\cor} {{\rm cor} \kern.15em}
\newcommand{\la}{\langle}
\newcommand{\ra}{\rangle}
\newcommand{\eps}{\varepsilon}

\newcommand{\xb}{\bar{x}}
\newcommand{\tos}{\rightrightarrows}

\begin{document}
\thispagestyle{empty}
\begin{center}
{\large\bf\sc Extended Forward-Backward Algorithm}
\end{center}

\begin{center}
  {\small\begin{tabular}{c}
  Marc LASSONDE and Ludovic NAGESSEUR\\
  LAMIA,
  Universit\'e des Antilles et de la Guyane\\
  97159 Pointe \`a Pitre, France\\
  E-mail: marc.lassonde@univ-ag.fr, ludovic.nagesseur@univ-ag.fr 
  \end{tabular}}
\end{center}

\medbreak\noindent
\textbf{Abstract.}
We propose an extended forward-backward algorithm for approximating a zero of
a maximal monotone operator which can be split as the extended sum of
two maximal monotone operators.
We establish the weak convergence in average of the sequence generated by
the algorithm under assumptions similar to those used in
classical forward-backward algorithms.
This provides as a special case an algorithm for solving
convex constrained minimization problems without qualification condition. 

\medbreak\noindent
\textbf{Keywords:} maximal monotone operator, $\eps$-enlargement, extended sum,
forward-backward algorithm, projected subgradient algorithm,
splitting algorithm, convergence in average.\bigbreak

\medbreak\noindent
\textbf{2010 Mathematics  Subject Classification:}
47J25, 47H05, 90C25, 

\section{Introduction}

Let $\Hcal$ be a real Hilbert space and let $T:\Hcal\tos \Hcal$ be a maximal monotone
operator. A fundamental problem is that of designing efficient algorithms
for constructing a zero of $T$, i.e.\ a solution of the inclusion
\begin{equation}\label{basic}
0\in Tx.
\end{equation}
Proximal point algorithms have been introduced by Martinet \cite{Mar70,Mar72}
to solve (\ref{basic}) in specific situations such as when
$T$ is the subdifferential $\partial f$ of a proper convex lower semicontinuous function $f$:
in this case the inclusion (\ref{basic}) is just the minimization of $f$ over $\Hcal$.
The general situation of an arbitrary maximal monotone operator $T$
was then considered by Rockafellar \cite{Roc76} and Brezis-Lions \cite{BL78}.
See, e.g., Peypouquet-Sorin \cite{PS10} for a survey on the topic.
 
When the maximal monotone operator $T$  admits a decomposition
of the form $T=A+B$ with $A$ and $B$ maximal monotone on $\mathcal{H}$,
it is more efficient to use splitting methods
which combine a step for $A$ and a step for $B$:
see, e.g., Combettes \cite{Com04}
for a unified presentation of these methods.
A prominent example is the forward-backward algorithm
considered by Passty \cite{Pas79},
which consists of a forward step on $B$ and a backward (proximal) step on $A$.
A special case of this method is the projected subgradient algorithm
aimed at solving constrained minimization problems.

There are important examples, however, where the problem to be solved
involves two maximal monotone operators $A$ and $B$ whose sum $A+B$
is not maximal monotone
so that the usual splitting methods do not apply. 
A simple instance of such a situation is the problem
\begin{equation}\label{un}
\mbox{Find } x\in C \mbox{ such that } f(x)=\min f(C)
\end{equation}
of minimizing a proper lower semicontinuous convex function $f:\Hcal\to \Rex$ over
a nonempty closed convex subset $C$ of $\Hcal$.
This problem can be rewritten as
\begin{equation}\label{deux}
\mbox{Find } x\in \Hcal \mbox{ such that } 0\in\partial (f+\delta_C)(x),
\end{equation}
where $\delta_C$ is the indicator of the set $C$.
When a qualification condition is satisfied, e.g. $0\in \cor (\dom f-C)$, then
$\partial (f+\delta_C)=\partial f+\partial\delta_C$, so (\ref{deux})
is equivalent to
\begin{equation}\label{trois}
\mbox{Find } x\in \Hcal \mbox{ such that } 0\in (\partial f+\partial\delta_C)(x),
\end{equation}
where $\partial f+\partial\delta_C$ is maximal monotone: in this case,
standard splitting methods apply.
But without a qualification condition, such a decomposition of
$\partial (f+\delta_C)$ as the sum
$\partial f+\partial\delta_C$ is not admissible. However, it is always true that
$\partial (f+\delta_C)$ can be split as the \textit{extended sum}
$\partial f\esum\partial\delta_C$ of the two maximal monotone operators
$\partial f$ and $\partial\delta_C$. Hence, (\ref{deux})
is always equivalent to
\begin{equation}\label{quatre}
\mbox{Find } x\in \Hcal \mbox{ such that } 0\in(\partial f\esum\partial\delta_C)(x),
\end{equation}
where $\partial f\esum\partial\delta_C$ is maximal monotone.

The aim of this paper is to present a modified forward-backward splitting algorithm to
solve the inclusion (\ref{basic}) for a maximal monotone operator $T$
which can be written as the extended sum $A\esum B$ of two maximal monotone operators.
This includes (\ref{un}) as a special case.
We establish the weak convergence in average of the sequence
generated by the algorithm to a zero of $T$
under hypotheses on the step sizes similar to those used in 
the classical splitting algorithms
in Bruck \cite{Bru77}, Lions \cite{Lio78} or Passty \cite{Pas79}
where the maximal monotone operator $T$
is equal to the usual sum $A+B$ of two maximal monotone operators.

Previous extensions of splitting algorithms 
to treat the case of the extended sum of
two maximal monotone operators rely heavily on the assumption that the
$\eps$-enlargement of one of the operators is bounded on the solution set
as in Moudafi-Th\'era \cite{MT01}.
Such a strong assumption is avoided here.

\section{Tools}

This section presents the main concepts and results we need in the sequel.
Most of these tools are valid in arbitrary Banach spaces, but we restrict ourselves to the
Hilbert space setting which suffices for our purpose.

Let $\Hcal$ be a real Hilbert space with inner product $\la .,.\ra$.
A set-valued operator $T:\Hcal\tos \Hcal$
is identified with its graph $T\subset \Hcal\times \Hcal$,
its domain is the set
$\D(T)=\{ x\in \Hcal : Tx\ne \emptyset\,\}$.
A set-valued operator $T:\Hcal\tos \Hcal$, or a subset $T\subset \Hcal\times \Hcal$,
is said to be {\it monotone} provided
\[
\langle v-u,y-x\rangle \geq 0,\quad  \forall (x,u),(y,v)\in T,
\]
and {\em maximal monotone} provided it is monotone and 
maximal (under set inclusion) in the family of all monotone sets contained
in $\Hcal\times \Hcal$.

Let $f:\Hcal\to \Rex$ be a convex lower semicontinuous function
which is {proper}, that is, its
{effective domain} $\dom f=\{x\in \Hcal: f(x)<+\infty\}$ is nonempty. 
For $\eps\ge 0$, the {\it $\eps$-subdifferential} $\partial _\eps f$
of $f$ at $x\in \Hcal$ is the set:
\[
\partial _\eps f(x):=\{u\in \Hcal: f(x)+\la u,y-x\ra\le f(y)+\eps,\ \forall y\in \Hcal\}.
\]
The case $\eps=0$ gives the {\it subdifferential} of $f$ at $x$,
denoted by $\partial f(x)$. For every $\eps>0$ and $x\in \dom f$,
the set $\partial_\eps f(x)$ is nonempty. On the
contrary, the set $\partial f(x)$ may be empty at some points $x\in \dom f$.
The subdifferential $\partial f$ is a basic example of a maximal monotone
operator on $\Hcal$ (see Moreau \cite{Mor65}).

The usual (Minkowski) sum $\partial f + \partial g$ of the
subdifferentials of two proper convex lower semicontinuous
functions $f,g:\Hcal\to \Rex$ is a monotone operator which
is not maximal monotone in general. Said differently,
$\partial f + \partial g$ is a monotone subset of $\partial(f+g)$
and the inclusion is proper in general without additional assumptions.
However, an exact formula for the subdifferential of the sum can be given in terms of
the sum of the $\eps$-subdifferentials of each function, namely
\begin{equation}\label{HP}
\partial (f+g)(x)=\bigcap_{{\varepsilon }>0}\overline{\partial _{\eps
}f(x)+\partial_{\eps}g(x)},\quad \forall x\in \Hcal.
\end{equation}
This formula is due to Hiriart-Urruty and Phelps \cite{HP93}.
It shows that the operator on the right end side is maximal monotone.

\smallbreak
Let now $T:\Hcal\tos \Hcal$ be an arbitrary monotone operator.
Given  $\eps\ge 0$, the {\em $\eps$-enlargement} of $T$ is the operator
$T^\eps:\Hcal\tos \Hcal$ defined as follows:
\begin{equation*}
T^\eps x:=\{u\in\Hcal: \langle v-u,y-x\rangle \ge
-\eps, \ \forall (y,v)\in T \}.
\end{equation*}
For a proper convex lower semicontinuous function $f$, one has $\partial _\eps f(x)\subset
(\partial f)^\eps(x)$.
For more details on this concept, see, e.g.,
\cite{BS99,GLR06,RT99,RT02} and the references therein.
The following inequality will be useful in the sequel (see \cite{BS99}):

\begin{proposition}\label{transportation}
Let $T:\Hcal\tos \Hcal$ be maximal monotone and let $\eps_1\ge 0$, $\eps_2\ge 0$. Then
\[
\langle v-u,y-x\rangle \ge -(\sqrt{\eps_1}+\sqrt{\eps_2})^2,\quad
\forall (x,u)\in T^{\eps_1}, (y,v)\in T^{\eps_2} 
\]
\end{proposition}

As for the subdifferentials, the usual sum $(A+B)(x)=Ax+Bx$, $x\in \Hcal$, 
of two maximal monotone operators $A,B$ defines a
monotone operator which is not maximal monotone in general.
Thus, motivated by the above formula (\ref{HP}), a concept of \textit{extended sum},
based on the notion of enlargement of a monotone operator, was proposed in
\cite{RT99,RT02}, namely:

\begin{equation}\label{RT}
(A\esum B)(x):= \bigcap_{\eps >0}\overline{A^\eps x+ B^\eps x},\quad \forall x\in \Hcal.
\end{equation}
The following result asserts that the extended sum $A\esum B$ is indeed
a monotone extension of the usual sum $A+ B$ (see \cite{GLR06}):

\begin{proposition}\label{teoextsumsub}
If $A,B:\Hcal\tos \Hcal$ are two maximal monotone operators,
the extended sum $A\esum B$ is a monotone operator containing the usual sum $A+B$.
\end{proposition}
It follows from Proposition \ref{teoextsumsub} and Formula (\ref{HP}) that
for any proper lower semicontinuous convex functions $f,g:\Hcal\to \Rex$
we have
\begin{equation}\label{RT2}
\partial(f+g)=\partial f\esum \partial g.
\end{equation}
This formula is due to Revalski and Th\'era \cite{RT02} (see also \cite{GLR06}).
It shows that the extended sum $\partial f\esum \partial g$ is a maximal monotone operator.

\section{Extended forward-backward algorithm}
In what follows, $\mathcal{H}$ is a real Hilbert space
with inner product $\la .,.\ra$ and unit ball $B_\Hcal$,
and $T:\mathcal{H}\tos \mathcal{H}$ is a maximal monotone operator
which can be split as the extended sum of two maximal monotone operators $A$ and $B$:
\begin{equation} \label{extsum}
Tx:=\underset{\eps>0}{\bigcap}\overline{A^{\eps}x+B^{\eps}x},\mbox{ for all } x\in\mathcal{H}.
\end{equation}
Our aim is to take advantage of this splitting to design
a forward-backward algorithm for constructing a zero of $T$,
i.e.\ a solution $x\in \mathcal{H}$ of the inclusion (\ref{basic}).
The set of solutions of (\ref{basic}) is supposed to be nonempty and is
denoted by $\Scal$.
\ps
Let $(\lambda_n)$ be a sequence of positive real numbers defining the step sizes of a discretization scheme, and let $\sigma_n:=\sum_{k=1}^n\lambda_k$ denote the step length
at stage $n$. As usual, $\sigma_n$ is supposed to grow to $+\infty$ as $n\to\infty$.
Given a sequence $(x_n)$ in $\mathcal{H}$
together with a sequence of step sizes $(\lambda_n)$, we
define the sequence $(\xb_n)$ of weighted averages by: 
\begin{equation}\label{formwa}
 \xb_n=\frac{1}{\sigma_n}\sum_{k=1}^n\lambda_k x_k.
\end{equation}
\ps
We study the \textit{extended forward-backward} algorithm given by the iteration:
\begin{equation}\label{fbalg}
\begin{cases}
& x_0\in \D(B^{\eps_0})\\
&    x_{n+1}=(I+\lambda_n A)^{-1}(x_{n}-\lambda_n u_n)
\quad \mbox{with } u_n\in B^{\eps_n}x_n,\quad \forall n\in\mathbb{N}.
\end{cases}
\end{equation}

Our main result is as follows:
\begin{theorem}\label{teoprincip}
Assume $\mathrm{D}(A)\subset \cap_{\eps>0}\mathrm{D}(B^\eps)$ and
let $(x_n)$ be the sequence generated by (\ref{fbalg}),
with $(\lambda_n)$ and $(\eps_n)$ any sequences of positive real numbers verifying
\begin{equation*}\label{param0}
\sum\lambda_n=\infty,\quad \sum\lambda_n^{4/3}<\infty,\quad \eps_n=\lambda_n^{1/3}.
\end{equation*}
Assume further that:
\ps\noindent
{\rm (H1)} The sequence $(\eps_n u_n)$ is bounded,

\noindent
{\rm (H2)} For every $x\in \Scal$ there is $M>0$ such that
$0\in\bigcap_{n}\,\overline{A^{\eps_n}x+B^{\eps_n}x\cap (M/\eps_n)B_\Hcal}$.

\ps\noindent
Then, the sequence $(\xb_n)$ of weighted averages given by (\ref{formwa})
weakly converges to a point in $\Scal$.
\end{theorem}

\proof
The argument relies on Passty's variant of Opial's lemma
(see \cite{Pas79,PS10}):
\medbreak\noindent
\textbf{Opial-Passty's Lemma}
{\em Let $(\lambda_n)$ be a sequence of positive numbers such that $\sum\lambda_n=\infty$, 
$(x_n)$ a sequence in a Hilbert space $\mathcal{H}$ and
$(\xb_n)$ its associated sequence of weighted averages given by (\ref{formwa}).
Let $\Scal\subset \mathcal{H}$ be a nonempty closed subset.
Assume

{\rm (i)} Every weak sequential cluster point of $(\xb_n)$ belongs to $\Scal$,

{\rm (ii)} $\lim \|x_n-x\|$ exists for every $x\in\Scal$.\\
Then, the sequence $(\xb_n)$ converges to a point in $\Scal$.}
\medbreak\noindent
The theorem is proved by showing that (i) and (ii) hold.
\ps
(i) Let $(x,y)\in T$, that is, $y\in\overline{A^{\eps}x+B^{\eps}x}$
for every $\eps>0$.
Fix $\eps>0$, and take $y_{\eps}=y_{1,\eps}+y_{2,\eps}$ with
$y_{1,\eps}\in A^{\eps}x$, $y_{2,\eps}\in B^{\eps}x$ and
$\|y-y_{\eps}\|\le \eps$.
As $u_n\in B^{\eps_n}x_n$ and $y_{2,\eps}\in B^{\eps}x$,
Proposition \ref{transportation} yields:
\begin{equation}\label{extsum2}
     \langle u_n-y_{2,\eps},x_n-x\rangle \geq -(\sqrt{\eps_n}+\sqrt{\eps})^2
     \geq -2(\eps_n+\eps)
\end{equation}
Since $x_n-\lambda_n u_n\in x_{n+1}+\lambda_n Ax_{n+1}$
by (\ref{fbalg}) and $y_{1,\eps}\in A^{\eps}x$, we get
\begin{equation*}
    \left\langle\frac{x_n-\lambda_n u_n-x_{n+1}}{\lambda_n}-y_{1,\eps}, x_{n+1}-x\right\rangle\geq -\eps,
\end{equation*}
equivalently:
\begin{equation}\label{3}
    2\langle x_n-x_{n+1},x_{n+1}-x\rangle\geq 2\lambda_n\langle u_n+y_{1,\eps},x_{n+1}-x\rangle-2\lambda_n\eps.
\end{equation}
By using the equality, valid for all $a, b, c\in\mathcal{H}$:
\begin{equation*}
    \|a-c\|^2=\|a-b\|^2+\|b-c\|^2+2\langle a-b,b-c\rangle,
\end{equation*}
we see that (\ref{3}) is equivalent to
\begin{equation}\label{3b}
 \|x_n-x\|^2-\|x_{n+1}-x\|^2\geq \|x_n-x_{n+1}\|^2
 +2\lambda_n\langle u_n+y_{1,\eps},x_{n+1}-x\rangle-2\lambda_n\eps.
\end{equation}
Now, combining the facts that:
\[
    2\lambda_n\langle u_n+y_{1,\eps},x_{n+1}-x_n\rangle 
    \geq -\|x_{n+1}-x_n\|^2-\lambda_n^2\|u_n+y_{1,\eps}\|^2,
\]
and
\begin{equation*}
    \langle u_n+y_{1,\eps},x_{n}-x\rangle=\langle u_n-y_{2,\eps},x_{n}-x\rangle+\langle y_\eps,x_{n}-x\rangle,
\end{equation*}
we derive from  (\ref{3b}) that:
\begin{equation}\label{3c}
\begin{array}{rcl}
    \|x_n-x\|^2-\|x_{n+1}-x\|^2
   & \geq &\|x_n-x_{n+1}\|^2-\|x_{n+1}-x_n\|^2-\lambda_n^2\|u_n+y_{1,\eps}\|^2\\
     && +2\lambda_n\langle u_n-y_{2,\eps},x_{n}-x\rangle+2\lambda_n\langle y_\eps,x_{n}-x\rangle-2\lambda_n\eps.
\end{array}
\end{equation}
Then, injecting (\ref{extsum2}) into (\ref{3c}) we obtain
\begin{equation}\label{ineqzer}
    2\lambda_n\langle y_\eps,x_{n}-x\rangle\leq\|x_n-x\|^2-\|x_{n+1}-x\|^2 
    +\lambda_n^2\|u_n+y_{1,\eps}\|^2+6\lambda_n\eps+4\lambda_n\eps_n.
\end{equation}
Summing up these inequalities for $n$ going from 1 to $k\in\mathbb{N}^*$,
and then dividing by $\sigma_k=\sum_{n=1}^k\lambda_n$, we get: 
\begin{equation}\label{final}
    2\langle y_\eps,\xb_k-x\rangle\leq \frac{\|x_1-x\|^2}{\sigma_k}
    +\frac{\sum_{n=1}^k\lambda_n^2\|u_n+y_{1,\eps}\|^2}{\sigma_k}
    +\frac{\sum_{n=1}^k\lambda_n(6\eps+4\eps_n)}{\sigma_k}.
\end{equation}
It follows from $\eps_n=\lambda_n^{1/3}$ and (H1) that
there exists a constant $M>0$ such that
$$
\lambda_n^2\|u_n\|^2=\lambda_n^{4/3}\eps_n^{2}\|u_n\|^2\le \lambda_n^{4/3} M^2,
$$
hence, since $\sum \lambda_n^{4/3}<\infty$,
$$
\lim_{k\to\infty}\,\sum_{n=1}^k\lambda_n^2\|u_n+y_{1,\eps}\|^2
\le 2\lim_{k\to\infty}\,\left (\sum_{n=1}^k\lambda_n^{4/3} M^2
+\sum_{n=1}^k\lambda_n^2\|y_{1,\eps}\|^2\right)<+\infty.
$$
Finally, since $\sigma_k\to +\infty$, we conclude that
\[
\lim_{k\to\infty}\,\frac{\sum_{n=1}^k\lambda_n^2\|u_n+y_{1,\eps}\|^2}{\sigma_k}=0.
\]
On the other hand, since $\eps_n\searrow 0$, there is $N\in\N$ such that
$\eps_n<\eps$ for $n>N$, hence, for $k$ large enough:
\[
\frac{\sum_{n=1}^k\lambda_n(6\eps+4\eps_n)}{\sigma_k}
\le \frac{\sum_{n=1}^N\lambda_n(6\eps+4\eps_n)}{\sigma_k}
+ 10\eps,
\]
from which we derive that
\[
 \limsup_{k\to\infty}\, \frac{\sum_{n=1}^k\lambda_n(6\eps+4\eps_n)}{\sigma_k}
\le 10\eps.
\]
Consequently, if $\overline{x}$ is any weak sequential cluster point of
the sequence $(\xb_n)$, letting $k\rightarrow\infty$ on both sides of (\ref{final}) yields:
\begin{equation}
    \langle y_\eps,\overline{x}-x\rangle\leq 5\eps,
\end{equation}
and then letting $\eps\rightarrow 0$ gives:
\begin{equation}\label{ineqzb}
    \langle y,x-\overline{x}\rangle\geq 0.
\end{equation}
Thus, every weak sequential cluster point $\bar{x}$ of $(\xb_n)$ satisfies 
(\ref{ineqzb}) for every $(x,y)\in T$.
The maximal monotonicity of $T$ then implies that
$0\in T\bar{x}$, that is, $\bar{x}\in\Scal$.
\ps
(ii)
Let $x\in\Scal$.
It follows from (H2) that there exist $M>0$ and sequences
$y_{1,n}\in A^{\eps_n}x$ and $y_{2,n}\in B^{\eps_n}x$
with $\eps_n\|y_{2,n}\|<M$ and $\|y_{1,n}+y_{2,n}\|<\eps_n/(1+\|x_n-x\|)$.
Without loss of generality, we may also assume that
$\eps_n\|y_{1,n}\|<M$ and, by (H1), that
$\eps_n\|u_n\|<M$ for all $n$.

Proceeding as above with $y=0$, $\eps=\eps_n$ and
$y_{\eps_n}=y_{1,n}+y_{2,n}$, we arrive at (\ref{ineqzer}),
hence
\begin{equation}\label{ineqzerf}
   2\lambda_n\langle y_{\eps_n},x_{n}-x\rangle
   \leq\|x_n-x\|^2- \|x_{n+1}-x\|^2 +2\lambda_n^{2}(\|u_n\|^2+
   \|y_{1,n}\|^2)+10\lambda_n\eps_n.
\end{equation}   
Then, using $\eps_n=\lambda_n^{1/3}$ and taking the above upper bounds into account, we get:
\begin{equation}\label{ineqzerfb}
   2\lambda_n\langle y_{\eps_n},x_{n}-x\rangle
   \leq\|x_n-x\|^2- \|x_{n+1}-x\|^2 +4\lambda_n^{4/3}M^2+10\lambda_n^{4/3},
\end{equation}
and, since $\|y_{\eps_n}\|\le \eps_n/(1+\|x_n-x\|)$,
\begin{equation}\label{ineqzerf2}
-2\lambda_n^{4/3}=-2\lambda_n\eps_n \leq -2\lambda_n\|y_{\eps_n}\|\|x_{n}-x\|
 \leq 2\lambda_n\langle y_{\eps_n},x_{n}-x\rangle,
\end{equation}
so, putting (\ref{ineqzerfb}) and (\ref{ineqzerf2}) together,
\begin{equation}\label{ineqzerf3}
\|x_{n+1}-x\|^2\le \|x_n-x\|^2+\lambda_n^{4/3}(4M^2+12).
\end{equation}
Since $\sum \lambda_n^{4/3}<\infty$, we derive from (\ref{ineqzerf3}) that
$\lim\|x_n-x\|^2$ exists hence also $\lim\|x_n-x\|$.
\cqfd

A useful special case of Theorem \ref{teoprincip} is as follows:
\begin{coro}\label{corprincip}
Assume $\mathrm{D}(A)\subset \cap_{\eps>0}\mathrm{D}(B^\eps)$ and
let $(x_n)$ be the sequence generated by (\ref{fbalg}),
with $(\lambda_n)$ and $(\eps_n)$ any sequences of positive real numbers verifying
\begin{equation*}\label{param}
\sum\lambda_n=\infty,\quad \sum\lambda_n^{4/3}<\infty,\quad \eps_n=\lambda_n^{1/3}.
\end{equation*}
Assume further that:
\ps\noindent
{\rm (H1)} The sequence $(\eps_n u_n)$ is bounded,

\noindent
{\rm (H2')} $\Scal=(A+B)^{-1}(0)$.

\ps\noindent
Then, the sequence $(\xb_n)$ of weighted averages given by (\ref{formwa})
weakly converges to a point in $\Scal$.
\end{coro}

\proof
We claim that Assumption (H2') implies Assumption (H2).
Indeed, let $x\in \Scal$. Then, we can find $v\in Ax\cap(-Bx)$, so
for every small $\eps_n>0$, 
$0\in Ax+Bx\cap (\|v\|/\eps_n)B_\Hcal$. Hence (H2) holds.
The result therefore follows from Theorem \ref{teoprincip}.
\cqfd

\begin{remark}\label{rem0}
{\rm
Examples of sequences $(\lambda_n)$ and $(\eps_n)$ verifying the assumptions
in the above results are $\lambda_n=1/n$ and $\eps_n=(1/n)^{1/3}$.
As was pointed out by the referee, the proof remains valid for any
sequences $(\lambda_n)$ and $(\eps_n)$ verifying the relations 
\[
\sum\lambda_n=\infty,\quad \sum(\lambda_n/\eps_n)^{2}<\infty,
\quad \sum\lambda_n\eps_n<\infty,\quad \eps_n\searrow 0.
\]
Clearly, the optimal combination 
for these relations to be satisfied is when
$(\lambda_n/\eps_n)^{2}=\lambda_n\eps_n$, that is, $\eps_n=\lambda_n^{1/3}$.}
\end{remark}

\begin{remark}\label{rem1}
{\rm
Assumption (H1) is of course satisfied if, as in \cite[Theorem 2]{Pas79},
the sequence $(u_n)$ is supposed to be bounded.
Simple examples show that (H1) may hold while
$(u_n)$ is not bounded, see Example \ref{exemple} below.
}
\end{remark}

\begin{remark}\label{rem2}
{\rm
Assumption (H2') is of course satisfied if, as in Passty's result \cite[Theorem 2]{Pas79},
the usual sum $A+B$ is maximal monotone. In the aforementioned theorem,
the assumptions are 
\[
\sum\lambda_n=\infty,\quad \sum\lambda_n^{2}<\infty,\quad u_n\in Bx_n,
\quad (u_n) \mbox{ bounded}.
\]
Passty's result and our results are therefore not comparable.
However, simple examples exist where (H2) is satisfied whereas (H2') is not,
see Example \ref{exemple} below. 
}
\end{remark}

\begin{remark}\label{rem3}
{\rm
Assumption (H2') is satisfied if, as in \cite[Theorem 2]{MT01}
(dealing with the double-backward scheme),
for every $x\in\Scal$ there exists $\eps>0$ such that $A^{\eps}x$ is bounded.
Indeed, let $x\in\Scal$. Since 
\[
0\in\bigcap_{n}\,\overline{A^{\eps_n}x+B^{\eps_n}x},
\]
there is a sequence $(y_{1,n},y_{2,n})\in A^{\eps_n}x\times B^{\eps_n}x$
such that $y_{1,n}+y_{2,n}\to 0$. From the boundedness of $A^{\eps}x$
and $A^{\eps_n}x\subset A^{\eps}x$ for large $n$,
we derive that $(y_{1,n})$ is bounded, hence we may assume that
$(y_{1,n})$ weakly converges to some point $y$ and $(y_{2,n})$ weakly converges to $-y$.
Since the sets $A^{\eps_n}x$ and $B^{\eps_n}x$ are weakly closed,
we derive that $y\in A^{\eps_n}x$ and $-y\in B^{\eps_n}x$ for every $n$.
From the maximal monotonicity of $A$ and $B$, we conclude that
$y\in Ax\cap (-Bx)$, that is, $x\in (A+B)^{-1}(0)$.
Thus, $\Scal=(A+B)^{-1}(0)$, hence, (H2') holds. 
}
\end{remark}
\begin{ex}\label{exemple}
{\rm
In $\Hcal=\R$, let
$f:\R\to\R$ be given by $f(x)=-\sqrt{x}$ for $x\ge 0$, $f(x)=+\infty$ otherwise,
and let $C=\{0\}$.
We consider the problem (\ref{un}) of minimizing
the proper lower semicontinuous convex function $f$
over the nonempty closed convex set $C$. As already seen in the introduction, this problem
amounts to the inclusion (\ref{basic}) where $T=\partial\delta_C\esum \partial f$
is maximal monotone by (\ref{RT2}).
Clearly, the set of solutions of this problem is the nonempty set $\Scal=\{0\}$.
Putting $A=\partial\delta_C$ and $B=\partial f$, the algorithm (\ref{fbalg}) becomes
the projected approximate subgradient algorithm, noting that $(I+\partial\delta_C)^{-1}$ is
the projection $P_C$ onto $C$:
\begin{equation}\label{fbalg2}
\begin{cases}
& x_0\in \dom f\\
&    x_{n+1}=P_C(x_{n}-\lambda_n u_n)
\quad \mbox{with}\quad u_n\in (\partial f)^{\eps_n}(x_n),\quad \forall n\in\mathbb{N}.
\end{cases}
\end{equation}
We show that the assumptions of Theorem \ref{teoprincip} are satisfied.

First, note that $\D(A)=\{0\}$ and $\D(\partial_\eps f)=\dom f=[0,+\infty[$, hence
$\D(A)\subset \bigcap_{\eps>0} \D(B^{\eps})$.

Next, we have to check the assumptions (H1) and (H2).
Easy computations show that $\partial_\eps \delta_C(0)={]}{-\infty}, +\infty[$
and $\partial_\eps f(0)={]}{-\infty}, -1/4\eps]$.
 
For (H1), observe that $x_n=0$ for
every $n=1,2,\ldots$, since $C=\{0\}$. We may therefore choose
$u_n= -1/4\eps_n\in \partial_{\eps_n} f(0)$ so that $\eps_nu_n=-1/4$.
Assumption (H1) is satisfied.

For (H2), observe that the only point in $\Scal$ is $x=0$. Since
$\partial_\eps f(0)\cap [-1/4\eps, 1/4\eps]=\{-1/4\eps\}$, we have
\[
0=1/4\eps-1/4\eps \in \partial_\eps \delta_C(0)+\partial_\eps f(0)\cap (1/4\eps)B_\R
\]
for every $\eps>0$. Assumption (H2) is satisfied.

Thus, we may apply Theorem \ref{teoprincip} to derive the convergence in average of the sequence generated by the algorithm (\ref{fbalg2}) to the point minimizing $f$ over $C$.
We point out that no sequence $u_n\in \partial_{\eps_n} f(x_n)$ is bounded,
the set $\partial f(x_n)=\partial f(0)$ is empty, hence
the operator $\partial \delta_C+\partial f$ is empty, a fortiori
the set $(\partial \delta_C+\partial f)^{-1}(0)$ is also empty,
the sets $\partial_\eps \delta_C(0)$ and $\partial_\eps f(0)$ are not bounded.
So none of the previous algorithms in \cite{Bru77,Lio78,Pas79,MT01} is applicable.
}
\end{ex}

{\small
\def\cdprime{$''$} \def\cprime{$'$}

}
\end{document}